\documentclass[preprints,article,accept,moreauthors,pdftex]{my-mdpi} 

\usepackage[labelformat=simple]{subfig}

\firstpage{1} 
\pubvolume{10}
\issuenum{4}
\articlenumber{290}
\pubyear{2021}
\copyrightyear{2021}
\externaleditor{Academic Editors: Stevan Pilipovi\'{c} and Chris Goodrich} 
\datereceived{1 August 2021} 
\daterevised{27 Aug and 22 Sept 2021} 
\dateaccepted{28 October 2021} 
\datepublished{1 November 2021} 
\hreflink{https://doi.org/10.3390/\newline axioms10040290} 
\doinum{10.3390/axioms10040290}

\pdfoutput=1

\usepackage{subfig}
\usepackage{mathtools,amssymb}
\usepackage{mathrsfs} 

\allowdisplaybreaks

\Title{Hybrid Method for Simulation of a Fractional COVID-19 Model with Real Case Application}

\TitleCitation{Hybrid Method for Simulation of a Fractional COVID-19 Model with Real Case Application}

\Author{Anwarud Din $^{1}$\orcidA{}, 
Amir Khan $^{2}$\orcidB{}, 
Anwar Zeb $^{3}$\orcidC{},
Moulay Rchid Sidi Ammi $^{4,}$*\orcidD{},
Mouhcine Tilioua $^{4}$\orcidE{}
and Delfim~F.~M.~Torres $^{5}$\orcidF{}}

\AuthorNames{Anwarud Din, Amir Khan, Anwar Zeb, Moulay Rchid Sidi Ammi, Mouhcine Tilioua and Delfim F. M. Torres}

\AuthorCitation{Din,~A.; Khan,~A.; Zeb,~A.; Sidi~Ammi,~M.R.; Tilioua,~M.; Torres,~D.F.M.}

\address{%
$^{1}$ \quad Department of Mathematics, Sun Yat-Sen University, 
Guangzhou 510275, China; anwarud@mail.sysu.edu.cn\\
$^{2}$ \quad Department of Mathematics and Statistics, 
University of Swat, Mingora 19130, 
Khyber Pakhtunkhwa, Pakistan; amirkhan@uswat.edu.pk\\
$^{3}$ \quad Department of Mathematics, Abbottabad Campus, 
COMSATS University Islamabad, Abbottabad 22060, 
Khyber~Pakhtunkhwa, Pakistan; anwar@cuiatd.edu.pk\\
$^{4}$ \quad MAIS Laboratory, AMNEA Group, FST Errachidia, 
Moulay Isma\"{\i}l University of Mekn\`es, \mbox{P.O. Box 509 Boutalamine,} 
Errachidia 52000, Morocco; m.tilioua@umi.ac.ma\\
$^{5}$ \quad Center for Research and Development in Mathematics and Applications (CIDMA), 
\mbox{Department of Mathematics,} University of Aveiro, 3810-193 Aveiro, Portugal; delfim@ua.pt}

\corres{Correspondence: sidiammi@ua.pt}

\abstract{In this research, we provide a mathematical analysis for the novel coronavirus
responsible for COVID-19, which continues to be a big source of threat for humanity. 
Our fractional-order analysis is carried out using a non-singular kernel type operator 
known as the Atangana--Baleanu--Caputo (ABC) derivative. We parametrize the model
adopting available information of the disease from Pakistan in the period 
9th April to 2nd June 2020. We obtain the required solution with the help of 
a hybrid method, which is a combination of the decomposition method 
and the Laplace transform. Furthermore, a sensitivity analysis 
is carried out to evaluate the parameters that are more sensitive 
to the basic reproduction number of the model. Our results are compared 
with the real data of Pakistan and numerical plots are presented 
at various fractional orders.}



\begin{document}


\section{Introduction}

The novel coronavirus SARS-CoV-2, responsible for COVID-19,
which is member of the family of Severe Acute Respiratory Syndrome (SARS) viruses,
has been recognized as the most dangerous virus of this decade~\cite{15g}. 
This virus has become the new novel strain of the SARS family, 
which was not recognized in humans before~\cite{15h}. COVID-19 has not just 
affected humans, but~also a number of animals have been infected by the virus. 
The SARS-CoV-2 virus has been transmitted from human to human and similarly in animals,
but its origin is still a controversy~\cite{Koopmans2021}. 
Infected humans and different species of various animals are 
recognized as active causes of spreading of the virus~\cite{15g}. 
In the past, some similar viruses, like the Middle East Respiratory Syndrome Coronavirus (MERS-CoV), 
were spread out from camels to human population, and~for SARS-CoV-1 the civet  
was recognized as the source of spreading into humans.
For COVID-19, the~main source or the major reason of spreading is
human-to-human interaction, where the virus transmission is easily
made by an infected person to a susceptible one. Currently, thousands of
research studies have been proposed and many predictions have been
given on COVID-19 dynamics, see~in \cite{MR4200529,LEMOSPAIAO2020100885,%
MR4173153,axioms10010018,axioms10030135} 
and references therein. Our paper is, however, different
from those in the literature. 
In~\cite{MR4200529}, special focus is given to the transmissibility 
of the so-called superspreaders, with~numerical simulations being 
given for data of Galicia, Spain, and~Portugal. It turns out
that, for~each region, the~order of the Caputo derivative takes 
a different value, always different from one, showing the relevance 
of considering fractional-order models to investigate COVID-19.
The work in~\cite{LEMOSPAIAO2020100885} studies the COVID-19 pandemic 
in Portugal until the end of the three states of emergency, describing
well what has happen in Portugal with respect to the evolution of
active infected and hospitalized individuals.
In~\cite{MR4173153,axioms10010018}, a~non-fractional but stochastic time-delayed model 
for COVID-19 is given, with~the aim to study
the situation of Morocco. In~\cite{axioms10030135}, the~authors provide
a $S$-$E$-$I$-$P$-$A$-$H$-$R$-$F$ model while here we propose
a much simpler $P$-$I$-$Q$ model (our model has only three state variables
while the model in~\cite{axioms10030135} is much more complex, 
with eight state variables). In~\cite{axioms10030135}, 
the authors use the classical operator of Caputo; 
differently, here we use the more recent 
ABC derivatives that, in~contrast, 
use non-singular kernels that allow 
us to consider a much simpler model.
While the main result in~\cite{axioms10030135} 
is the proof of the global stability of the 
disease free equilibrium; in contrast, 
here we prove Ulam--Hyers stability.
We also construct a practical algorithm 
to compute numerically the solution of the model
(see Section~\ref{sec:5}), while such algorithmic approach is not 
addressed in~\cite{axioms10030135}. 
Moreover, we do a sensitivity analysis
to the parameters of the model.
Such sensitivity analysis is also not addressed 
in~\cite{axioms10030135}. 
In contrast with the work in~\cite{axioms10030135}, that investigates
the realities of Wuhan, Spain, and Portugal,
we study the case of~Pakistan.

COVID-19 generally transfers by interaction with humans 
in close contact for a particular time period, 
with most common symptoms of sneezing and coughing. 
The virus droplets stay on the layer of
matters and when they come to the contact with any susceptible
human, then the virus symptoms easily transfer to the individuals.
Such infected humans can pass the infection to others by touching their mouth,
eyes, or~nose. This virus has the strength to be alive on different
surfaces, like cardboard and copper, for~many hours up to some days. 
As the time passes, the~amount of the virus symptoms decreases over a
time span and might not be alive in sufficient amount for spreading
the infection. It has been recorded that the symptoms appearance
and COVID-19 infection initial stage lies between 1 to 14~days~\cite{15g}. 
Several countries have prepared and implemented 
a COVID-19 vaccine program and are trying to protect their populations. 
However, to~date, there is yet no treatment available. 
At present, the~most effect way to protect ourselves from the virus 
remains the quarantine or isolation, effective use of mask, following 
the guidelines that have been passed by governments of all countries along with 
the World Health Organization (WHO).

Modeling of infectious diseases has a rich literature and a number
of research articles have been developed, both using classical dynamical systems 
as well as fractional models~\cite{MR4200529,axioms10030135}.
Fractional-order derivatives can be useful and helpful as compared to classical derivatives,
because the dynamics of real phenomena can be comprehensively
understood by fractional-order derivatives due to its special
properties, i.e.,~hereditary and memory~\cite{1,2,6,7,8,9,10,11}. 
For a comparison between classical (integer-order)
and fractional-order models, see in~\cite{MR4200529,axioms10030135}.
Roughly speaking, ordinary derivatives cannot distinguish the phenomenon 
at two distinct closed points. To~sort out this problem of ordinary derivatives, 
generalized derivatives have been introduced in the framework of fractional calculus~\cite{3}. 
The first concept of fractional-order derivative was given by Leibniz and L'H\^opital in
1695. Aiming quantitative analysis, optimization, and~numerical
estimations, many number of attempts have been made by employing
fractional differential equations (FDEs)
\cite{12,13,14,15,15a,15b,15c,15d,15e,15f,15g,15h,15i,15j,15k,15l,15m,15n,15o,15p,15q,15r,15s}.
The growing interest in the modeling of complex real-world issues
with the use of FDEs is due to its numerous properties that can not
be found in the ordinary sense. These characteristics allow FDEs 
to model effectively not only non-Markovian processes but also 
non-Gaussian phenomena~\cite{math9161883}. 
Different non-classical fractional-order derivatives and
different kinds of FDEs were proposed~\cite{4,5,MR4260534}. 
Among them, one has the Atangana--Baleanu--Caputo (ABC) derivative,
which is a nonlocal fractional derivative with a non-singular kernel,
connected with various applications. For~a discussion of the ABC 
and related operators see in~\cite{16,MR3891804}, and for~their use on contemporary 
modeling we refer the interested reader to~the works in \cite{17,18,19,ammi2021necessary}.

The famous method of decomposition was developed from the 1970s 
to the 1990s by George Adomian, to~analytically handle nonlinear problems. 
After that, the~Adomian decomposition method became a powerful tool to simulate
analytical or approximate solutions for various problems of an
applied nature. Many mathematical models have been studied by the
applications of homotopy, Laplace Adomian Decomposition Method (LADM), 
and variational methods~\cite{20,21,22}. To~the best of our knowledge, 
no one has studied a variable order epidemic model with ABC 
derivatives by the LADM. Motivated by this fact, here we study 
a fractional-order COVID-19 epidemic model with ABC derivatives by the
Laplace Adomian decomposition algorithm. In~particular, we use Banach and
Krassnoselskii fixed point theorems to define some sufficient
conditions to prove existence and uniqueness of solution.
As stability is important for the estimated solution, we 
consider Ulam type stability through nonlinear functional
analysis. The~aforementioned stability is investigated for
ordinary fractional derivatives in many research papers, see, e.g.,~ in
\cite{24,25,26}, but~research on Ulam type stability regarding 
ABC derivatives is a rarity. At~the end of the paper, our results 
are illustrated with real data based on Pakistan COVID-19 cases 
in March~2020. 

The paper is organized as follows. Section~\ref{sec:2} 
is devoted to the model formulation. Section~\ref{sec:3} is concerned
with some preliminary results on fractional differential equations.
Existence and uniqueness are carried out in Section~\ref{sec:4}. 
Section~\ref{sec:5} deals with the solution of the COVID-19 model using 
the LADM. Some plots are given in Section~\ref{sec:6}, showing the simplicity 
and reliability of the proposed algorithm. In~Section~\ref{sec:7}, 
a sensitivity analysis is given to find the most sensitive parameter
with respect to the basic reproduction number. We end with Section~\ref{sec:8} 
of conclusions, including some possible future directions of~research.


\section{Model~Formulation}
\label{sec:2}

Mathematical modeling plays a major role in investigating and thus
controlling the dynamics of a disease, particularly in the
vaccination privation or at the initial phases of the epidemic.
Several mathematical models can be found in~\cite{7,8,9,10}. We
formulated a fractional COVID-19 epidemic model, similar to other
diseases~\cite{20,21}, and~predict its future behavior. Inspired by
FDEs using the ABC derivative, we aim to simulate
the COVID-19 transmission in the form of
\begin{equation}
\left\{
\begin{split}
\label{m1}
^{\mathcal{ABC}}\mathbf{D}_{t}^{\theta}{P}(t)
&=\lambda-\gamma
P(t)I(t)-d_0P(t),\\
^{\mathcal{ABC}}\mathbf{D}_{t}^{\theta}{I}(t)
&=\gamma P(t)I(t)-(d_0+h+\eta)I(t)+\sigma Q(t),\\
^{\mathcal{ABC}}\mathbf{D}_{t}^{\theta}{Q}(t)
&=\eta I(t)-(d_0+\mu+\sigma) Q(t),
\end{split}\right.
\end{equation}
along with initial conditions
\begin{equation}
P(0)=P_0,\ I(0)=I_0,\  Q(0)=Q_0,
\end{equation}
where $^{\mathcal{ABC}}\mathbf{D}_{t}^{\theta}$ is the ABC fractional
derivative of order $0<\theta\leq 1$ (see Definition~\ref{d1} 
in Section~\ref{sec:3}). In~this model, $P(t)$ represents the amount
of susceptible humans, $I(t)$ stands for the population of infected
humans, and~$Q(t)$ represents the population of quarantined humans
at time $t$. The meaning of the parameters of model \eqref{m1}
are given in Table~\ref{table:01}. We take the below assumptions to the given~system:
\begin{itemize}
\item[$A_1$.] All the variables and parameters of the system are non-negative.
\item[$A_2$.] The susceptible people transfer to the infectious compartment 
with a constant susceptible inflow into population.
\item[$A_3$.] Originally infectious or susceptible persons transfer 
to the quarantined class while reported cases return to the infected class 
from quarantined classes.
\end{itemize}

The basic reproduction number $R_0$, which represents the secondary cases
for the model \eqref{m1}, is easily demonstrated to be given by
\begin{eqnarray}
\label{eq2}
R_0=\frac{\gamma \lambda
(d_0+\mu+\sigma)}{d_0(d_0+\mu+\sigma)(d_0+h)+\eta(d_0+\mu)}.
\end{eqnarray}
\begin{specialtable}[H]
\setlength{\tabcolsep}{9.7mm}
\caption{Parameters description defined in the given model \eqref{m1}.
\label{table:01}}
\begin{tabular}{ll}
\toprule
\textbf{Notation}  & \textbf{Description} \\\midrule
$\lambda$ & Rate of recruitment \\
$\gamma$ & Transmission rate of disease \\
$d_0$ & Natural death rate  \\
$\eta$ & Transmission rate of infected to quarantine \\
$\mu$ & Deaths in quarantined zone \\
$\sigma$ &Transmission flow of quarantined to become infectious  \\
$h$ & Rate of deaths in infected zone\\\bottomrule
\end{tabular}
\end{specialtable}

In addition, $I(t)+P(t)+Q(t)=\mathrm{N}(t),$ where $N$ represents
the total~population.


\section{Preliminary~Results}
\label{sec:3}

For completeness, here we recall necessary definitions 
and results from the~literature.

\begin{Definition}[See~\cite{ammi2021necessary,6}]
\label{d1} 
If $x$ is an absolutely continuous function
and $0<\theta\leq{1}$, then the ABC derivative is given by
\begin{eqnarray}
\label{d2}
^{\mathcal{ABC}}\mathbf{D}_{t}^\theta\phi(t)
=\frac{\mathcal{ABC}(\theta)}{1-\theta}
\int_0^t\frac{d}{dy}x(\omega)\mathcal{M}_\theta\left[
\frac{-\theta}{1-\theta}\left(t-\omega\right)^\theta\right],
\end{eqnarray}
where $\mathcal{ABC}(\theta)$ is a normalization function such that 
$1=\mathcal{ABC}(1)=\mathcal{ABC}(0)$ and
$\mathcal{M}_\theta$ is a special Mittag--Leffler function.
\end{Definition}

\begin{Remark}
By replacing $\mathcal{M}_\theta\bigg[\frac{-\theta}{1-\theta}\bigg(t-\omega\bigg)^\theta\bigg]$ 
with $\mathcal{M}_1=\exp\bigg[\frac{-\theta}{1-\theta}\bigg(t-\omega\bigg)\bigg]$ one obtains 
the so-called Caputo--Fabrizio derivative. Additionally, we have
$$
^\mathcal{ABC}\mathbf{D}_{0}^\theta[constant]=0.
$$ 
\end{Remark}

\begin{Remark}
\label{dl} 
Let $x(t)$ be a function  having fractional ABC derivative.
Then, the Laplace transform of $^{\mathcal{ABC}}\mathbf{D}_{0}^\theta x(t)$ 
is given by
\begin{equation*}
\mathscr{L}\left[^{\mathcal{ABC}}\mathbf{D}_{0}^\theta x(t)\right]
=\frac{\mathcal{ABC}(\theta)}{[s^\theta(1-\theta)+\theta]}\left[s^\theta
\mathscr{L}[x(t)]-s^{\theta-1}x(0)\right].
\end{equation*}
\end{Remark}

\begin{Lemma}[See~\cite{24}]
\label{w11}  
The solution to 
\begin{gather*}
^{\mathcal{ABC}}\mathbf{D}_{0}^\theta x(t)
=z(t),\quad t\in [0, T],\\
x(0)=x_0,
\end{gather*} 
$1>\theta>0$, is given by
$$
x(t)=x_0+\frac{(1-\theta)}{\mathcal{ABC}(\theta)}z(t)
+\frac{\theta}{\mathcal{ABC}(\theta)\Gamma(\theta)}
\int_0^t (t-\omega)^{\theta-1}z(\omega)d\omega.
$$
\end{Lemma}

\begin{Theorem}[See~\cite{26}]
\label{p9}
Let $\mathbf{X}=C[0, T]$ and consider the Banach space defined by
$\mathbf{Z}=\mathbf{X}\times \mathbf{X}\times \mathbf{X}$ with 
the norm-function $\|A\|=\|(P,I,Q)\|
=\max_{t\in [0, T]}[|P(t)+|I(t)|+|Q(t)|]$. Consider $\mathbf{B}$ 
to be a convex subset of $\mathbf{Z} $ and   
$\mathbf{F}$, and $\mathbf{G}$ be operators such~that
\begin{enumerate}
\item $\mathbf{F}u + \mathbf{G}u\in \mathbf{B}$  $\forall$ $ u\in \mathbf{B}$,
\item $\mathbf{F}$ is a contraction, and
\item $\mathbf{G}$ is compact and continuous.
\end{enumerate} 

Then, $\mathbf{F}u+\mathbf{G}u=u$ possesses at least one solution.
\end{Theorem}


\section{Qualitative Analysis of the Proposed~Model}
\label{sec:4}

Here, we rewrite the right-hand sides of \eqref{m1} as
\begin{equation}
\label{Ntt}
\begin{split}
\mathrm{f}_1(t,P(t), I(t), Q(t))
&=-\gamma I(t)P(t)+\lambda-d_0P(t),\\
\mathrm{f}_2(t,P(t), I(t), Q(t))
&=\gamma I(t)P(t)-(d_0+h+\eta)I(t)+\sigma Q(t), \\
\mathrm{f}_3(t,P(t), I(t), Q(t))
&=\eta I(t)-(d_0+\sigma+\mu) Q(t).
\end{split}
\end{equation}

By using \eqref{Ntt}, we have
\begin{equation} 
\label{r2}
\begin{split}
^{\mathcal{ABC}}\mathbf{D}_{+0}^\theta \mathcal{A}(t)
&=\Phi(t, \mathcal{A}(t)),
\quad t\in [0, \tau], 
\quad 0<\theta\leq 1,\\
\mathcal{A}(0)&=\mathcal{A}_0.
\end{split}
\end{equation}

In view of Lemma~\ref{w11}, \eqref{r2} yields
\begin{multline}
\label{21}
\mathcal{A}(t)=\mathcal{A}_0(t)
+\left[\Phi(t, \mathcal{A}(t))-\Phi_0(t)\right]
\frac{(1-\theta)}{\mathcal{ABC}(\theta)}\\
+\frac{\theta}{\Gamma(\theta) \mathcal{ABC}(\theta)}
\int_0^t (t-\omega)^{\theta-1} \Phi(\omega, \mathcal{A}(\omega))d\omega,
\end{multline}
where
\begin{equation}
\label{200}
\begin{gathered}
\mathcal{A}(t)=\left\{
\begin{split}
&P(t)\\&
I(t)\\&
Q(t)
\end{split}
\right., \ \ \
\mathcal{A}_0(t)=\left\{
\begin{split}
&P_0\\&
I_0\\&
Q_0
\end{split}
\right.,\\
\Phi(t, \mathcal{A}(t))=\left\{
\begin{split}
&\mathrm{f}_1(t,P, I, Q)\\&
\mathrm{f}_2(t,P, I, Q)\\&
\mathrm{f}_3(t,P, I, Q).
\end{split}
\right.,\ \ \
\Phi_0(t)=\left\{
\begin{split}
&\mathrm{f}_1(0,P_0, I_0, Q_0)\\&
\mathrm{f}_2(0,P_0, I_0, Q_0)\\&
\mathrm{f}_3(0,P_0, I_0, Q_0).
\end{split}
\right.
\end{gathered}
\end{equation}

Using \eqref{21} and \eqref{200}, 
we define the two operators $\mathbf{F}$ and $\mathbf{G}$  
as follows:
\begin{equation}
\label{23}
\begin{split}
\mathbf{F}(\mathcal{A})
&= \mathcal{A}_0(t)+\left[\Phi(t, \mathcal{A}(t))-\Phi_0(t)\right]
\frac{(1-\theta)}{\mathcal{ABC}(\theta)},\\
\mathbf{G}(\mathcal{A})&= 
\frac{\theta}{\Gamma(\theta) \mathcal{ABC}(\theta)}
\int_0^t (t-\omega)^{\theta-1}\Phi(\omega, \mathcal{A}(\omega))d\omega.
\end{split}
\end{equation}

For existence and uniqueness, we assume some basic axioms 
and a Lipschitz hypothesis:
\begin{itemize}
\item[(H1)] there is $C_{\Phi}$ and $D_\Phi$ such that
$$ 
|\Phi(t, \mathcal{A}(t))|\leq C_{\Phi}\|\mathcal{A}\|+D_{\Phi};
$$
\item[(H2)] there is $L_{\Phi}>0$ such that $\forall$ 
$\mathcal{A},\ \mathcal{\bar{A}}\in \mathbf{Z}$ one has
$$
|\Phi(t, \mathcal{A})-\Phi(t, \mathcal{\bar{A}})|
\leq L_{\Phi}[\|\mathcal{A}\|-\|\mathcal{\bar{A}}\|].
$$
\end{itemize}

\begin{Theorem}
Under hypotheses $(H1)$ and $(H_2)$, Equation \eqref{21} possesses at least
one solution, which implies that \eqref{m1} possesses an equal number of solutions 
if $\frac{(1-\theta)}{\mathcal{ABC}(\theta)}L_{\Phi}<1$.
\end{Theorem}

\begin{proof} 
The theorem is proved in two steps, with~the help of Theorem~\ref{p9}.
(i) Consider  $\mathcal{\bar{A}}\in \mathbf{B}$, 
where $\mathbf{B}=\{\mathcal{A}\in \mathbf{Z}:\ \|\mathcal{A}\|\leq \rho, \rho>0\}$ 
is a closed and convex set. Then, for~$\mathbf{F}$ in \eqref{23}, we have
\vspace{12pt}\begin{equation}
\label{p1}
\begin{split}
\|\mathbf{F}(\mathcal{A})-\mathbf{F}(\mathcal{\bar{A}})\| 
&= \frac{(1-\theta)}{\mathcal{ABC}(\theta)}\max_{t\in [0, \tau]}\left|
\Phi(t, \mathcal{A}(t))-\Phi(t, \mathcal{\bar{A}}(t))\right|\\
&\leq \frac{(1-\theta)}{\mathcal{ABC}(\theta)}L_{\Phi}\|\mathcal{A}-\mathcal{\bar{A}}\|.
\end{split}
\end{equation}

Therefore, $\mathbf{F}$ is a contraction.
(ii) We want $\mathbf{G}$ to be relatively compact. 
For that it suffices that $\mathbf{G}$ is equicontinuous and bounded.
Obviously, $\mathbf{G}$ is continuous as $\Phi$ is continuous and 
for all $\mathcal{A}\in \mathbf{B}$ one has
\begin{equation}
\label{q12}
\begin{split}
\|\mathbf{G}(\mathcal{A})\|
&=\max_{t\in [0, \tau]}\left|\frac{\theta}{\Gamma(\theta)\mathcal{ABC}(\theta)}
\int_0^t (t-\omega)^{\theta-1}\Phi(\omega, \mathcal{A}(\omega))d\omega\right|\\
&\leq \frac{\theta}{\Gamma(\theta)\mathcal{ABC}(\theta)}
\int_0^{\tau}(\tau-\omega)^{\theta-1}|\Phi(\omega, \mathcal{A}(\omega))|d\omega\\
&\leq\frac{\tau^\theta}{\mathcal{ABC}(\theta)\Gamma(\theta)}[C_{\Phi}\rho+D_{\Phi}].
\end{split}
\end{equation}

Thus, \eqref{q12} shows the boundedness of $\mathbf{G}$. 
For equi-continuity, we assume $t_1>t_2 \in [0, \tau]$, so that
\begin{equation}
\label{n120}
\begin{split}
|&\mathbf{G}(\mathcal{A}(t_1))-\mathbf{G}(\mathcal{A}(t_2))|\\
&=\frac{\theta}{\mathcal{ABC}(\theta)\Gamma(\theta)}\left|
\int_{0}^{t_1}(t_1-\omega)^{\theta-1}\Phi(\omega, \mathcal{A}(\omega))d\omega
-\int_{0}^{t_2}(t_2-\omega)^{\theta-1}\Phi(\omega, \mathcal{A}(\omega))d\omega\right|\\
&\leq\frac{[C_{\Phi}\rho+D_{\Phi}]}{\mathcal{ABC}(\theta)\Gamma(\theta)}|t_1^\theta-t_2^\theta|.
\end{split}
\end{equation}

The right-hand side in \eqref{n120} goes to zero at $t_1\rightarrow t_2$. 
Remembering that $\mathbf{G}$ is continuous, 
$$ 
|\mathbf{G}(\mathcal{A}(t_1))-\mathbf{G}(\mathcal{A}(t_2))|\rightarrow 0
\ \text{ as }\  t_1\rightarrow t_2.
$$

Having the boundedness and continuity of $\mathbf{G}$, we conclude that
$\mathbf{G}$ is uniformly continuous and bounded. According to the
theorem of Arzel\'{a}--Ascoli, $\mathbf{G}$ is relatively
compact and therefore entirely continuous. It follows from Theorem~\ref{p9}
that the integral Equation \eqref{21} has at least one solution.
\end{proof}

Now, we show~uniqueness.

\begin{Theorem}
\label{90} 
Under hypotheses $(H1)$ and $(H_2)$, Equation \eqref{21} possesses a
unique solution and this implies that \eqref{m1} possesses also 
a unique solution if $\frac{(1-\theta)L_\Phi}{\mathcal{ABC}(\theta)}
+\frac{\tau^\theta L_\Phi}{\mathcal{ABC}(\theta)\Gamma(\theta)}<1$.
\end{Theorem}

\begin{proof}
Let the operator $\mathbf{T}:\mathbf{Z}\rightarrow \mathbf{Z}$
be defined by
\begin{multline}
\label{w12}
\mathbf{T}\mathcal{A}(t)=\mathcal{A}_0(t)+\bigg[\Phi(t,
\mathcal{A}(t))-\Phi_0(t)\bigg]\frac{(1-\theta)}{\mathcal{ABC}(\theta)}\\
+\frac{\theta}{\mathcal{ABC}(\theta)\Gamma(\theta)}\int_0^t
(t-\omega)^{\theta-1} \Phi(\omega, \mathcal{A}(\omega))d\omega,
\quad t\in [0, \tau].
\end{multline}
Let $\mathcal{A}, \mathcal{\bar{A}} \in \mathbf{Z}$. 
Then, one can take
\begin{equation}
\label{pio}  
\begin{split}
\|\mathbf{T}\mathcal{A}-\mathbf{T}\mathcal{\bar{A}}\|
&\leq\frac{(1-\theta)}{\mathcal{ABC}(\theta)}\max_{t\in [0, \tau]}\left|
\Phi(t, \mathcal{A}(t))-\Phi(t, \mathcal{\bar{A}}(t))\right|\\
&\quad +\frac{\theta}{\Gamma(\theta)\mathcal{ABC}(\theta)}
\max_{t\in [0, \tau]}\bigg|\int_0^t (t-\omega)^{\theta-1} 
\Phi(\omega, \mathcal{A}(\omega))d\omega\\
&\qquad -\int_0^t (t-\omega)^{\theta-1} 
\Phi(\omega, \mathcal{\bar{A}}(\omega))d\omega\bigg|\\
&\leq\Xi\|\mathcal{A}-\mathcal{\bar{A}}\|,
\end{split}
\end{equation}
where
\begin{equation}
\label{u126}
\Xi=\frac{(1-\theta)L_\Phi}{\mathcal{ABC}(\theta)}
+\frac{\tau^\theta L_\Phi}{\Gamma(\theta)
\mathcal{ABC}(\theta)}.
\end{equation}

Thus, $\mathbf{T}$ is a contraction from \eqref{pio}. Therefore,
\eqref{21} possesses a unique solution.
\end{proof}

Next, in~order to investigate the stability of our problem,
we consider a small disturbance $\phi \in C[0, T]$, 
with $\phi(0)=0$, that depends only on the~solution.

\begin{Lemma}
\label{lemma3.3}
Let $\phi \in C[0, T]$  with $\phi(0)=0$ such that
$|\phi(t)|\leq \varepsilon$ for $\varepsilon >0$
and consider the problem
\begin{equation}
\label{w2}
\begin{split}
^{\mathcal{ABC}}\mathbf{D}_{+0}^\theta \mathcal{A}(t)
&=\Phi(t, \mathcal{A}(t))+\phi(t),\\
\mathcal{A}(0)&=\mathcal{A}_0.
\end{split}
\end{equation}

The solution of \eqref{w2} satisfies the following relation:
\begin{equation}
\label{u8}
\begin{split}
\bigg|\mathcal{A}(t)
&-\bigg(\mathcal{A}_0(t)+\bigg[
\Phi(t, \mathcal{A}(t))-\Phi_0(t)\bigg]\frac{(1-\theta)}{\mathcal{ABC}(\theta)}\\
&\quad +\frac{\theta}{\mathcal{ABC}(\theta)\Gamma(\theta)}
\int_0^t (t-\omega)^{\theta-1} \Phi(\omega, \mathcal{A}(\omega))d\omega\bigg)\bigg|\\
&\leq \frac{\Gamma(\theta)+ \tau^\theta}{\Gamma(\theta)\mathcal{ABC}(\theta)}\varepsilon
=\Omega_{\tau,\theta}.
\end{split}
\end{equation}
\end{Lemma}

\begin{proof}
The proof is standard and is omitted here.
\end{proof}

\begin{Theorem}
Consider hypotheses $(H1)$ and $(H_2)$ 
along with \eqref{u8} of Lemma~\ref{lemma3.3}. 
Then, the solution to Equation \eqref{21} is Ulam--Hyers stable 
if $ \Xi <1$, where $\Xi$ is defined by \eqref{u126}.
\end{Theorem}

\begin{proof}
Assume $\mathcal{A}\in \mathbf{Z}$ 
and let $\mathcal{\bar{A}}\in \mathbf{Z}$ be the
unique solution of \eqref{21}. Then,
\begin{eqnarray}
\|\mathcal{A}-\mathcal{\bar{A}}\|
&=&\max_{t\in [0, T]}\bigg| \mathcal{A}(t)-\bigg(\mathcal{A}_0(t)
+\bigg[\Phi(t, \mathcal{\bar{A}}(t))
-\Phi_0(t)\bigg]\frac{(1-\theta)}{\mathcal{ABC}(\theta)}\nonumber\\
&\quad& +\frac{\theta}{\Gamma(\theta)\mathcal{ABC}(\theta)}
\int_0^t (t-\omega)^{\theta-1} \Phi(\omega,
\mathcal{\bar{A}}(\omega))d\omega\bigg)\bigg|\nonumber\\
&\leq&\max_{t\in [0, T]}\bigg| \mathcal{A}(t)-\bigg(\mathcal{A}_0(t)
+\bigg[\Phi(t, \mathcal{{A}}(t))
-\Phi_0(t)\bigg]\frac{(1-\theta)}{\mathcal{ABC}(\theta)}\nonumber\\
&\quad& +\frac{\theta}{\Gamma(\theta)\mathcal{ABC}(\theta)}
\int_0^t (t-\omega)^{\theta-1} \Phi(\omega, 
\mathcal{{A}}(\omega))d\omega\bigg)\bigg|\nonumber\\
&\quad& +\bigg|\bigg(\mathcal{A}_0(t)+\bigg[\Phi(t, \mathcal{{A}}(t))
-\Phi_0(t)\bigg]\frac{(1-\theta)}{\mathcal{ABC}(\theta)}\nonumber\\
&\quad& +\frac{\theta}{\Gamma(\theta)\mathcal{ABC}(\theta)}
\int_0^t (t-\omega)^{\theta-1} \Phi(\omega, \mathcal{{A}}(\omega))d\omega\bigg)\nonumber\\
&\quad& -\bigg(\mathcal{A}_0(t)+\bigg[\Phi(t, \mathcal{\bar{A}}(t))
-\Phi_0(t)\bigg]\frac{(1-\theta)}{\mathcal{ABC}(\theta)}\nonumber\\
&\quad& +\frac{\theta}{\Gamma(\theta)\mathcal{ABC}(\theta)}
\int_0^t (t-\omega)^{\theta-1} \Phi(\omega, \mathcal{\bar{A}}(\omega))d\omega\bigg)\bigg|\nonumber\\
&\leq&\Omega_{\tau, \theta}+\frac{(1-\theta)L_\Phi}{\mathcal{ABC}(\theta)}\|\mathcal{A}
-\mathcal{\bar{A}}\|+\frac{\tau^\theta L_\Phi}{\Gamma(\theta)
\mathcal{ABC}(\theta)}\|\mathcal{A}-\mathcal{\bar{A}}\|\nonumber\\
&\leq&\Omega_{\tau, \theta}+\Xi\|\mathcal{A}-\mathcal{\bar{A}}\|. \label{d23}
\end{eqnarray}

From \eqref{d23}, we can write that
\begin{equation}
\label{d45}
\|\mathcal{A}-\mathcal{\bar{A}}\|
\leq \frac{\Omega_{\tau,\theta}}{1-\Xi}\left\|
\mathcal{A}-\mathcal{\bar{A}}\right\|.
\end{equation}

The proof is complete.
\end{proof}


\section{Construction of an Algorithm for Deriving the Solution of the~Model}
\label{sec:5}

Herein, we derive a general series-type solution for the 
proposed system with ABC derivatives. Taking the Laplace transform in
model \eqref{m1}, we transform both sides of each equation and
we use the initial conditions to obtain that
\begin{equation}
\label{a3}
\left\{
\begin{split}
\mathscr{L}[P(t)]&=\frac{P_0}{s}+\frac{[s^\theta(1-\theta)
+\theta]}{s^\theta\mathcal{ABC}(\theta)}\mathscr{L}\left[
\lambda-\gamma P(t)I(t)-d_0P(t)\right],\\
\mathscr{L}[I(t)]&=\frac{I_0}{s}+\frac{[s^\theta(1-\theta)
+\theta]}{s^\theta\mathcal{ABC}(\theta)}
\mathscr{L}\left[\gamma I(t)P(t)-I(t)(d_0+h+\eta)+\sigma Q(t)\right],\\
\mathscr{L}[Q(t)]&=\frac{Q_0}{s}+\frac{[s^\theta(1-\theta)
+\theta]}{s^\theta\mathcal{ABC}(\theta)}\mathscr{L}\left[
\eta I(t)-(d_0+\mu+\sigma) Q(t)\right].
\end{split}\right.
\end{equation}

Now, considering each solution in the form of series,
\begin{equation}
\label{q1}
P(t)=\sum_{n=0}^{\infty}P_n(t),
\quad I(t)=\sum_{n=0}^{\infty}I_n(t),
\quad Q(t)=\sum_{n=0}^{\infty}Q_n(t),
\end{equation}
we separate the nonlinear term $P(t)I(t)$  
in terms of Adomian polynomials 
as
\begin{equation}
\label{q2}
P(t)I(t)=\sum_{n=0}^{\infty} H_n(t),\ \text{ where }\ 
H_n(t)=\frac{1}{n!}\frac{d^n}{d\lambda^n}\left[
\sum_{k=0}^{n}\lambda^kP_k(t)\sum_{k=0}^{n}
\lambda^kI_k(t)\right]\bigg|_{\lambda=0}.
\end{equation}

Therefore, from~\eqref{q1} and \eqref{q2},  
we obtain from \eqref{a3} that
\begin{equation}
\left\{
\begin{split}
\label{a30}
&\mathscr{L}\left[\sum_{n=0}^{\infty}P_n(t)\right]
=\frac{P_0}{s}+\frac{[s^\theta(1-\theta)
+\theta]}{s^\theta\mathcal{ABC}(\theta)}
\mathscr{L}\left[\lambda-\gamma \sum_{n=0}^{\infty}H_n(t)
-d_0\sum_{n=0}^{\infty}P_n(t)\right],\\
&\mathscr{L}\left[\sum_{n=0}^{\infty}I_n(t)\right]=\frac{I_0}{s}\\
&\qquad +\frac{[s^\theta(1-\theta)+\theta]}{s^\theta\mathcal{ABC}(\theta)}
\mathscr{L}\left[\gamma \sum_{n=0}^{\infty}H_n(t)-(d_0+h+\eta)
\sum_{n=0}^{\infty}I_n(t)+\sigma \sum_{n=0}^{\infty}Q_n(t)\right],\\
&\mathscr{L}\left[\sum_{n=0}^{\infty}Q_n(t)\right]
=\frac{Q_0}{s}+\frac{[s^\theta(1-\theta)
+\theta]}{s^\theta\mathcal{ABC}(\theta)}
\mathscr{L}\left[\eta \sum_{n=0}^{\infty}I_n(t)
-(d_0+\mu+\sigma) \sum_{n=0}^{\infty}Q_n(t)\right].
\end{split}\right.
\end{equation}

Now, comparing the terms on both sides of \eqref{a30}, one has
\begin{equation}
\label{a300}
\left\{
\begin{split}
\mathscr{L}[P_0(t)]&=\frac{P_0}{s},
\quad \mathscr{L}[I_0(t)]=\frac{I_0}{s},
\quad \mathscr{L}[Q_0(t)]=\frac{Q_0}{s}, \\
\mathscr{L}[P_1(t)]&=\frac{[s^\theta(1-\theta)+\theta]}{s^\theta
\mathcal{ABC}(\theta)}\mathscr{L}\left[\lambda-\gamma H_0(t)-d_0P_0(t)\right],\\
\mathscr{L}[I_1(t)]
&=\frac{[s^\theta(1-\theta)+\theta]}{s^\theta\mathcal{ABC}(\theta)}
\mathscr{L}\left[\gamma H_0(t)-(d_0+h+\eta)I_0(t)+\sigma Q_0(t)\right],\\
\mathscr{L}[Q_1(t)]
&=\frac{[s^\theta(1-\theta)+\theta]}{s^\theta\mathcal{ABC}(\theta)}
\mathscr{L}\left[\eta I_0(t)-(d_0+\mu+\sigma) Q_0(t)\right],\\
\mathscr{L}[P_2(t)]&=\frac{[s^\theta(1-\theta)+\theta]}{s^\theta
\mathcal{ABC}(\theta)}\mathscr{L}\bigg[\lambda
-\gamma H_1(t)-d_0P_1(t)\bigg],\\
\mathscr{L}[I_2(t)] &=\frac{[s^\theta(1-\theta)+\theta]}{s^\theta\mathcal{ABC}(\theta)}
\mathscr{L}\left[\gamma H_1(t)-(d_0+h+\eta)I_1(t)+\sigma Q_1(t)\right],\\
\mathscr{L}[Q_2(t)]
&=\frac{[s^\theta(1-\theta)+\theta]}{s^\theta\mathcal{ABC}(\theta)}
\mathscr{L}\left[\eta I_1(t)-(d_0+\mu+\sigma) Q_1(t)\right],\\
& \vdots\\
\mathscr{L}[P_{n+1}(t)]
&=\frac{[s^\theta(1-\theta)+\theta]}{s^\theta
\mathcal{ABC}(\theta)}\mathscr{L}\left[\lambda
-\gamma H_n(t)-d_0P_n(t)\right],\ n\geq 0,\\
\mathscr{L}[I_{n+1}(t)] &=\frac{[s^\theta(1-\theta)+\theta]}{s^\theta
\mathcal{ABC}(\theta)}\mathscr{L}\left[\gamma
H_n(t)-(d_0+h+\eta)I_n(t)+\sigma Q_n(t)\right],\ n\geq 0,\\
\mathscr{L}\left[Q_{n+1}(t)\right]
&=\frac{[s^\theta(1-\theta)+\theta]}{s^\theta\mathcal{ABC}(\theta)}
\mathscr{L}\left[\eta
I_n(t)-(d_0+\mu+\sigma) Q_n(t)\right],
\quad n\geq 0.
\end{split}\right.
\end{equation}

Applying the inverse Laplace transform to \eqref{a300}, we obtain that
\vspace{10pt}\end{paracol}
\nointerlineskip
\begin{equation}
\label{a301}
\left\{
\begin{split}
P_0(t)&=P_0,\quad I_0(t)=I_0,\quad Q_0(t)=Q_0,\\
P_1(t)&=\mathscr{L}^{-1}\left[\frac{[s^\theta(1-\theta)
+\theta]}{s^\theta\mathcal{ABC}(\theta)}\mathscr{L}\left[
\lambda-\gamma H_0(t)-d_0P_0(t)\right]\right],\\
I_1(t) &=\mathscr{L}^{-1}\left[\frac{[s^\theta(1-\theta)+\theta]}{s^\theta
\mathcal{ABC}(\theta)}\mathscr{L}\left[\gamma
H_0(t)-(d_0+h+\eta)I_0(t)+\sigma Q_0(t)\right]\right],\\
Q_1(t)&=\mathscr{L}^{-1}\left[\frac{[s^\theta(1-\theta)+\theta]}{s^\theta
\mathcal{ABC}(\theta)}\mathscr{L}\left[\eta
I_0(t)-(d_0+\mu+\sigma) Q_0(t)\right]\right],\\
P_2(t)&=\mathscr{L}^{-1}\left[
\frac{[s^\theta(1-\theta)+\theta]}{s^\theta\mathcal{ABC}(\theta)}
\mathscr{L}\left[\lambda-\gamma
H_1(t)-d_0P_1(t)\right]\right],\\ 
I_2(t)&=\mathscr{L}^{-1}\left[
\frac{[s^\theta(1-\theta)+\theta]}{s^\theta\mathcal{ABC}(\theta)}
\mathscr{L}\left[\gamma H_1(t)-(d_0+h+\eta)I_1(t)+\sigma Q_1(t)\right]\right],\\
Q_2(t) &=\mathscr{L}^{-1}\left[\frac{[s^\theta(1-\theta)
+\theta]}{s^r\mathcal{ABC}(\theta)}\mathscr{L}\left[\eta
I_1(t)-(d_0+\mu+\sigma) Q_1(t)\right]\right],\\
& \vdots\\
P_{n+1}(t)&=\mathscr{L}^{-1}\left[
\frac{[s^\theta(1-\theta)+\theta]}{s^\theta\mathcal{ABC}(\theta)}\mathscr{L}\left[
\lambda-\gamma H_n(t)-d_0P_n(t)\right]\right],\quad n\geq 0,\\
I_{n+1}(t) &=\mathscr{L}^{-1}\left[
\frac{[s^\theta(1-\theta)+\theta]}{s^\theta\mathcal{ABC}(\theta)}
\mathscr{L}\left[\gamma H_n(t)-(d_0+h+\eta)I_n(t)
+\sigma Q_n(t)\right]\right],\quad n\geq 0,\\
Q_{n+1}(t) &=\mathscr{L}^{-1}\left[
\frac{[s^\theta(1-\theta)+\theta]}{s^\theta
\mathcal{ABC}(\theta)}\mathscr{L}\left[\eta
I_n(t)-(d_0+\mu+\sigma) Q_n(t)\right]\right],\quad n\geq 0.
\end{split}\right.
\end{equation}
\begin{paracol}{2}
\switchcolumn


\section{Numerical Interpretation and~Discussion}
\label{sec:6}

To illustrate the dynamical structure of our infectious disease
model, we now consider a practical case study under various
numerical observations and given parameter values. The~concrete
parameter values we have used are shown in Table~\ref{tab2}. 
\begin{specialtable}[H]
\setlength{\tabcolsep}{4.1mm}
\caption{Numerical values for the parameters of model \eqref{m1}.}
\label{tab2}
\begin{tabular}{lll}
\toprule
\textbf{Notation}  & \textbf{Parameters Description} &\textbf{Numerical Value} \\
\midrule
$\lambda$ & Rate of recruitment & $0.003$\\ 
$\gamma$ & Transmission rate of disease & $0.009$\\ 
$d_0$ & Natural death rate & $0.009$\\ 
$\eta$ & Transmission rate of infected to quarantine& $0.004$\\ 
$\mu$ & Death rate in quarantine &$0.004$ \\
$\sigma$ & Transmission flow of quarantined to infectious & $0.003$ \\ 
$h$ & Rate of death for infected & $0.007$\\ 
$P_0$ & Initial population of susceptible& $10$ millions \\ 
$I_0$ & Initially infected population& $0.01$ millions \\ 
$Q_0$ & Quarantined population at $t=0$ & $0.0011$ millions\\
\bottomrule
\end{tabular}
\end{specialtable}

We assume that the initial susceptible, infected,  
and isolated populations are 10, 0.01, and 0.0011 million, respectively. 
Among the 21,000 selected population, the~density of susceptible population is
about 0.6 percent, the~infected population is 0.2~percent, and~the
isolated population is 0.2~percent.

By using the parameter values in Table~\ref{tab2}, we computed the first
three terms of the general series solution \eqref{a301} with
$\mathcal{ABC}(\theta)=1$ as
\begin{equation}
\label{ser123}
\begin{gathered}
\begin{split}
P(t) &=0.6+1.99712\bigg[1-\theta+\frac{\theta t^\theta}{\Gamma(\theta)}\bigg]\\
&\quad -0.00681\bigg[(1-\theta)^2t
+\frac{\theta^2 t^{2\theta}}{\Gamma(2\theta+1)}
+\frac{2\theta(1-\theta) t^{\theta+1}}{\Gamma(\theta+2)} \bigg]+\cdots,
\end{split}\\
\begin{split}
I(t) &= 0.2-0.5976\left[1-\theta+\frac{ \theta t^\theta}{\Gamma(\theta)}\right]\\
&\quad -0.003096\left[(1-\theta)^2t
+\frac{2\theta(1-\theta) t^\theta}{\Gamma(\theta)}
+\frac{\theta^2 t^{2\theta}}{\Gamma(2\theta+1)}\right]+\cdots,
\end{split}\\
\begin{split}
Q(t) &= 0.2+0.14\left[1-\theta+\frac{\theta t^\theta}{\Gamma(\theta)}\right]\\
&\quad -0.004058\left[(1-\theta)^2t+\frac{2\theta(1-\theta) 
t^\theta}{\Gamma(\theta)}
+\frac{\theta^2 t^{2\theta}}{\Gamma(2\theta+1)}\right]+\cdots.
\end{split}
\end{gathered}
\end{equation}

We have utilized the numeric computing environment MATLAB, 
version 2016, and~plotted the solution \eqref{a301} in
Figure~\ref{fig.03} by considering the first 
fifteen terms of the series \eqref{q1}.

Figure~\ref{fig.03} shows the dynamics of each one of the state variables 
in the classical sense when $\theta = 1$ (green curves). Similarly, 
by considering the model in ABC sense, that is, for~$\theta \in (0,1)$,
we plot in Figure~\ref{fig.03} each of the state variables to analyze 
the changes in comparison with the classical case. From~Figure~\ref{fig.03a},
we see that as we increase the order $\theta$ of the fractional ABC derivative, 
the susceptibility increases. Further, all fractional order derivatives
shows no effect after about 60 days, i.e.,~the susceptible
population stabilizes. Figure~\ref{fig.03b} shows that 
infected individuals tend to increase, with~different rates, when we
decrease the fractional-order: the smaller the fractional order $\theta$,
faster the increase rate, and vice~versa. All obtained curves for infected 
individuals, for~different values of the fractional order derivatives, 
approach towards a non-zero steady state, which shows that the disease 
will persist in the community if not properly managed. On~the other hand, 
Figure~\ref{fig.03c} shows that during the first month the disease 
progress with more and more people getting quarantined, irrespective 
of the order of the derivative. However, after~that, the~quarantined 
population tends to decline and, at~the end, there will be no quarantined 
individuals in the~community.

In Section~\ref{sec6.1}, we show that the
fractional model \eqref{m1} with ABC derivatives has the ability to describe 
effectively the dynamics of transmission of the current COVID-19~outbreak.
\begin{figure}[H]
{\captionsetup{position=bottom,justification=centering}
\subfloat[$P(t)$---susceptible individuals along time $t$]{
\label{fig.03a}\includegraphics[scale=0.45]{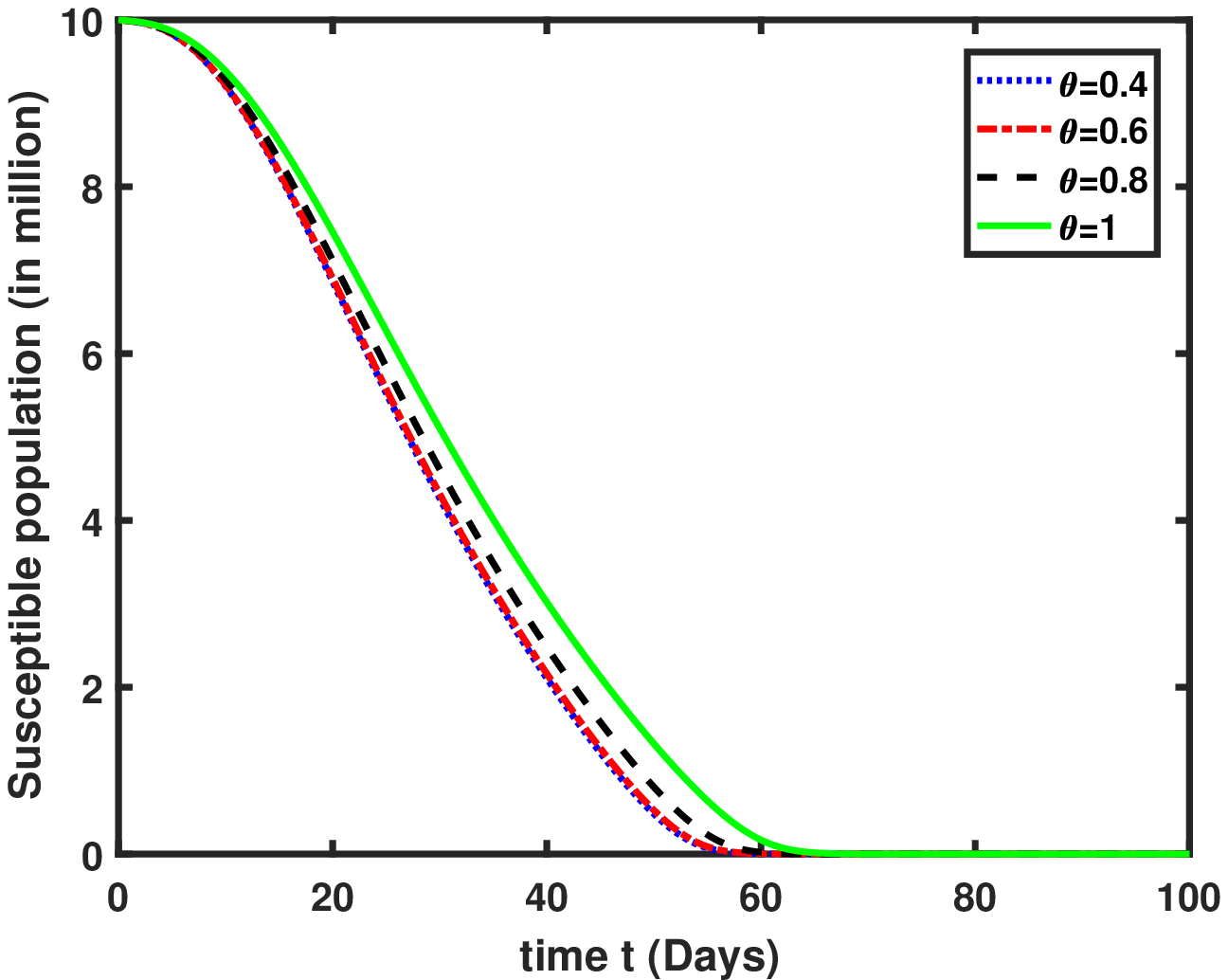}}
\subfloat[$I(t)$---infected individuals along time $t$]{
\label{fig.03b}\includegraphics[scale=0.45]{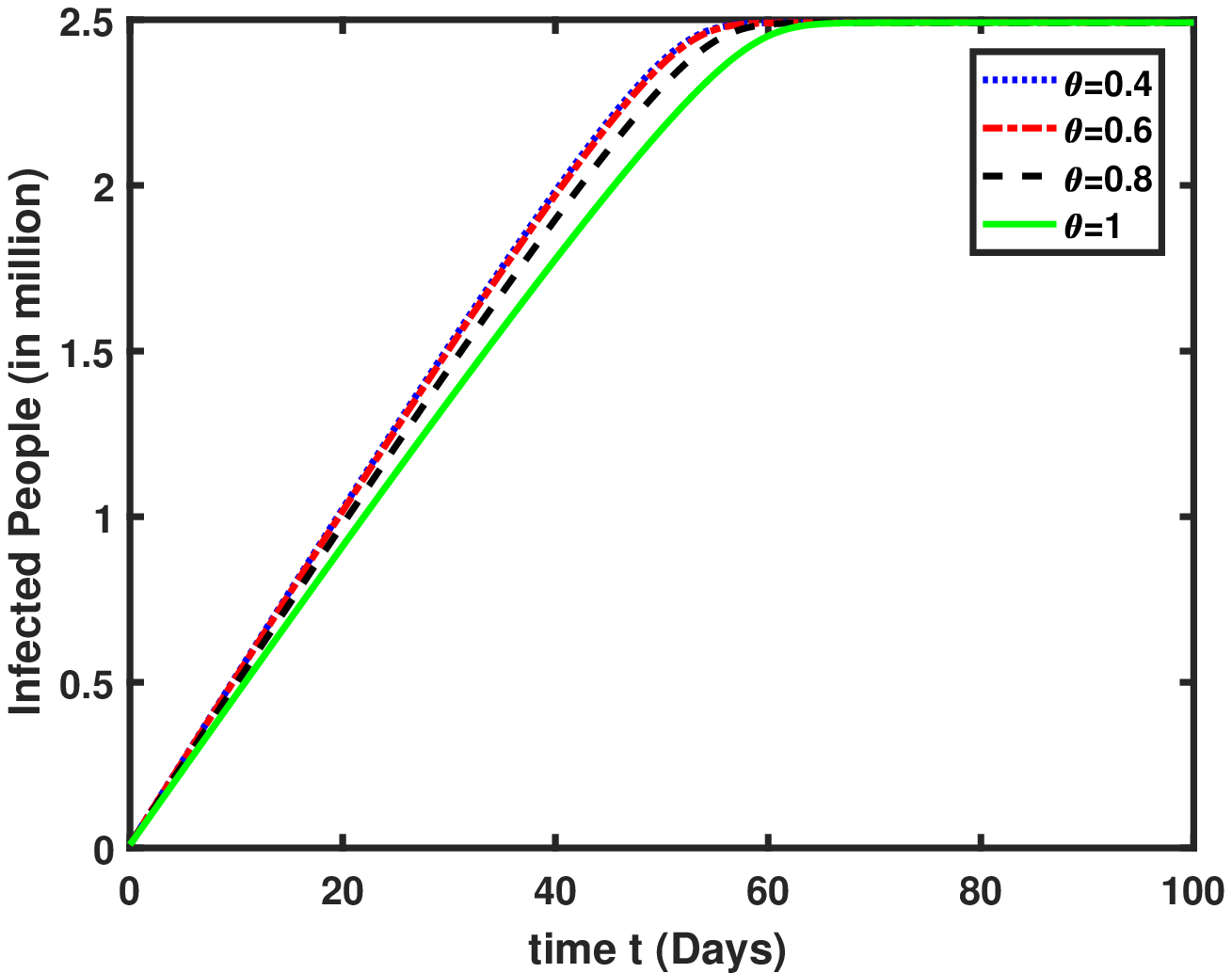}}\\
\subfloat[$Q(t)$---quarantined individuals along time $t$]{
\label{fig.03c}\includegraphics[scale=0.45]{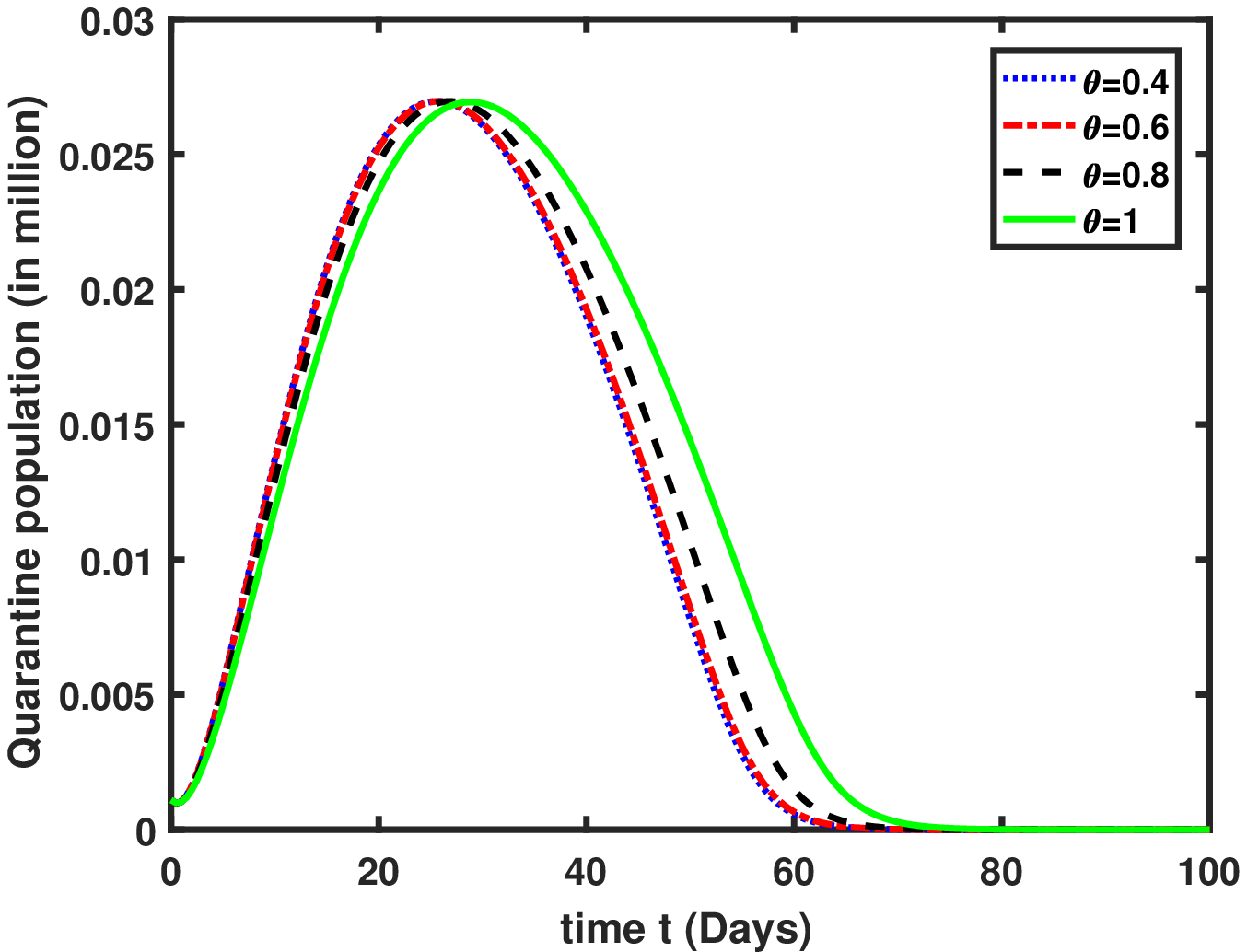}}}
\caption{Dynamical nature of susceptible, infected 
and quarantined individuals of the fractional ABC model \eqref{m1} 
for different values of the fractional-order $\theta$.}
\label{fig.03}
\end{figure}


\subsection{Case Study with Real Data: Khyber Pakhtunkhawa (Pakistan)}
\label{sec6.1}

The Khyber Pakhtunkhawa Province,
like other provinces of Pakistan and the rest of the world, 
is also being affected by COVID-19. 
We decided to calibrate our model with real data
of COVID-19 from Khyber Pakhtunkhawa, Pakistan,
from 9th April to the 2nd of June 2020. For~that,
we have used the minimization method of MATLAB
taking the initial weights 
\begin{equation*}
P(0) = 35,525,047,\quad 
I(0) = 10,485,\quad 
Q(0) = 18000,
\end{equation*}
determined from~the work in \cite{web}, and~$\theta=1$,
from which we arrived to the values of the parameters shown
in Table~\ref{Table3}.

Figure~\ref{fig5} shows the total number of individuals
infected by COVID-19 as registered from 9th April to the 2nd in
June 2020, which corresponds to the period of one month and 24 days 
used to calibrate our model.

Figure~\ref{fig6} compares the actual/real data of COVID-19 
with the curve of infected given by model \eqref{m1},
clearly showing the appropriateness of our model
to describe the COVID-19 outbreak.
\begin{specialtable}[H]
\setlength{\tabcolsep}{16.2mm}
\caption{{Parameter values for the case of Khyber Pakhtunkhawa, Pakistan.}}
\label{Table3}
\begin{tabular}{lll}
\toprule
\textbf{Notation}  & \textbf{Value} & \textbf{Reference} \\
\midrule
$\lambda$  & 0.028  & \cite{web} \\ 
$\gamma$ & {0.2}& Estimated \\ 
$d_0$  &  0.011 & \cite{web} \\ 
$\mu$  &0.2 & Estimated \\ 
$h$  & 0.06 & \cite{web} \\ 
$\sigma$  & 0.04 & Estimated  \\ 
$\eta$  & 0.3 & Estimated  \\  \bottomrule
\end{tabular}
\end{specialtable}
\vspace{-12pt}
\begin{figure}[H]
\includegraphics[width=3.5in]{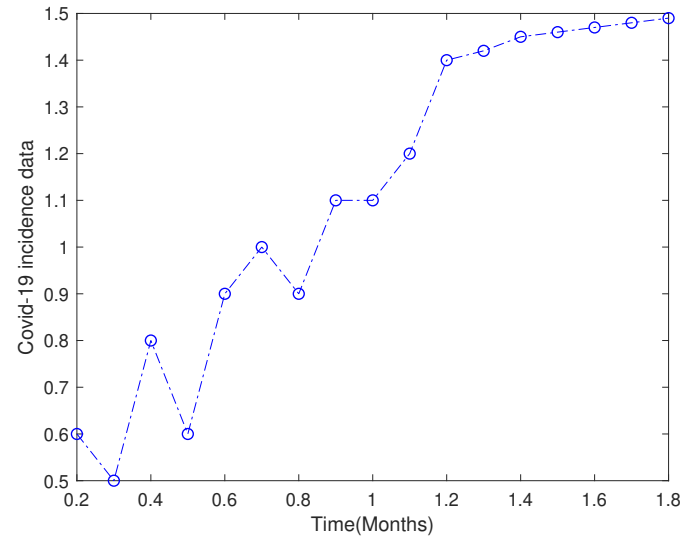}
\caption{Real data of infected individuals by COVID-19 
from Khyber Pakhtunkhwa, Pakistan, from~9th April 
to the 2nd of June 2020.}
\label{fig5}
\end{figure}
\vspace{-12pt}
\begin{figure}[H]
\includegraphics[width=3.5in]{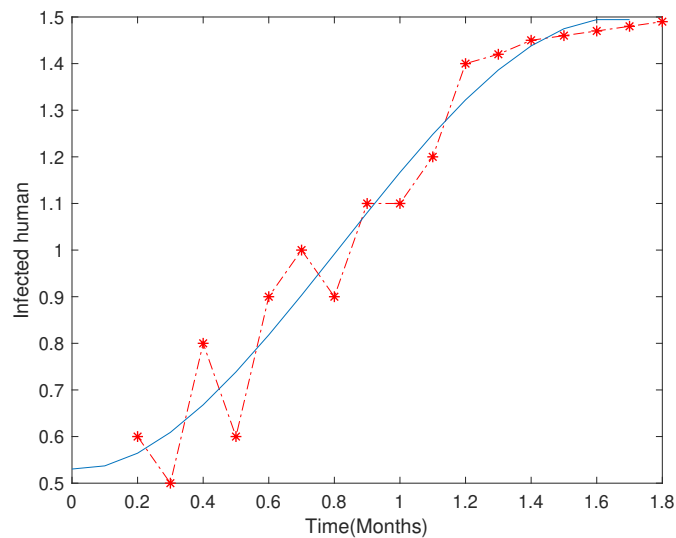}
\caption{Comparison of infected individuals by COVID-19: 
model \eqref{m1} output (in blue) versus real data 
of Khyber Pakhtunkhawa, Pakistan, from 9th April 
to the 2nd of June 2020 (in red).}
\label{fig6}
\end{figure}

Figure~\ref{fig7} projects the long-term behavior
of the COVID-19 outbreak during a period of eight months. 
We can see the data matches during the first 1.8 months  
and, additionally, we observe that the long-term behavior 
consists on a rise of infected individuals with time.
This means that if the government did not apply
proper strategies, the~incidence could increase 
drastically in the coming months.
\begin{figure}[H]
\includegraphics[width=3.5in]{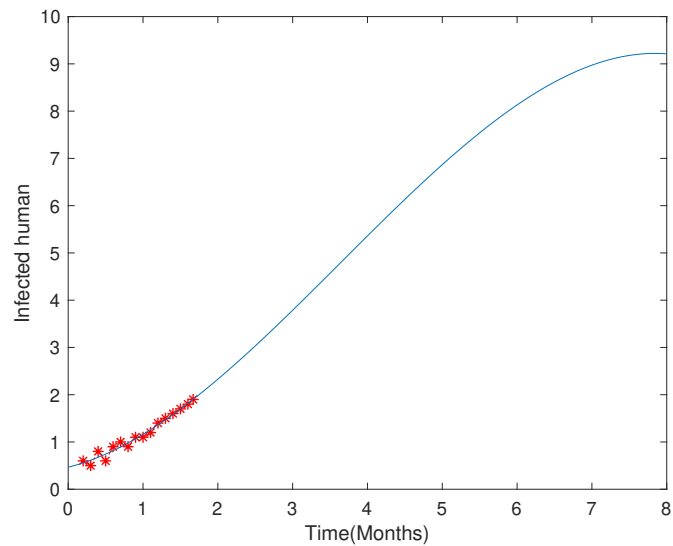}
\caption{Real data of infected individuals 
by COVID-19 in Khyber Pakhtunkhwa, Pakistan (first 1.8~months, in~red) 
and prediction from model \eqref{m1} during a period of 8 months (in blue).}
\label{fig7}
\end{figure}


\section{Sensitivity~Analysis}
\label{sec:7}

Here, we conduct a sensitivity analysis to evaluate the parameters 
that are sensitive in minimizing the propagation of the ailment. 
Although its computation is tedious for complex biological models, 
forward sensitivity analysis is recorded as an important component 
of epidemic modeling: the ecologist and epidemiologist gain a lot 
of insight from the sensitivity study of the basic reproduction number 
$R_0$ \cite{MR3557143}. In~Definition~\ref{def:sensitI}, we assume 
that the basic reproduction number $R_0$ 
is differentiable with respect to parameter $\omega$. Given \eqref{eq2},
this means that Definition~\ref{def:sensitI} makes sense for 
$\omega \in \{\gamma, \lambda, d_0, \mu, \sigma, h, \eta\}$.

\begin{Definition}
\label{def:sensitI}
The normalized forward sensitivity index of $R_0$
with respect to parameter $\omega$ is defined by
\begin{equation}
\label{eq:S:omega}
S_\omega=\frac{\omega}{R_0}\frac{\partial R_0}{\partial \omega}.
\end{equation}
\end{Definition}

As we have an analytical form for the basic reproduction number, recall \eqref{eq2},
we apply the direct differentiation process given in \eqref{eq:S:omega}. 
Not only do the sensitivity indexes show us the impact of various factors 
associated with the spread of the infectious disease, but they also provide us 
with valuable details on the comparative change between $R_{0}$ and the parameters.
Moreover, they also assist in the production of control strategies~\cite{MR3771538}.

Table~\ref{table4} demonstrates that $\gamma$, $h$, and $\sigma$ 
parameters have a positive effect on the basic reproduction number $R_0$, 
which means that the growth or decay of these parameters by 10\% 
would increase or decrease the reproduction
number by 10\%, 6.36\%, and 0.31\%, respectively. On~the other
hand, $d_0$-, $\mu$-, and $\eta$-sensitive indexes indicate that
increasing their values by 10\% would decrease the basic
reproduction number $R_0$ by 14.89\%, 0.09\%, and 1.68\%,
respectively.
\begin{specialtable}[H]
\setlength{\tabcolsep}{2.4mm}
\caption{Sensitivity indexes of the basic reproduction 
number $R_0$ \eqref{eq2} (see Definition~\ref{def:sensitI})
for relevant parameters of model \eqref{m1}.}
\label{table4}
\begin{tabular}{cccccc} \toprule
\textbf{Parameters} & \textbf{Sensitivity} & \textbf{Value} 
& \textbf{Parameters} & \textbf{Sensitivity} & \textbf{Value}  \\ \midrule
$\gamma$ & $S_\gamma$ & 1.00000000 &  $h$ & $S_{h}$     &  0.63636363 \\ \midrule
$d_0$  & $S_{d_0}$  &  $-$1.48944805 &  $\mu$ & $S_{\mu}$ & $-$0.00974026 \\ \midrule
$\sigma$ & $S_{\sigma}$  & 0.03165584 & $\eta$  & $S_{\eta}$  & $-$0.16883117 \\ \bottomrule
\end{tabular}
\end{specialtable}

The sensitivity of the basic reproduction number $R_0$
is also seen graphically in \mbox{Figure~\ref{fig:sensitiv}}.
\vspace{-12pt}
\begin{figure}[H]
{\captionsetup{position=bottom,justification=centering}
\subfloat[$R_{0}$ versus $\gamma$ and $d$]{
\label{fig8}\includegraphics[scale=0.45]{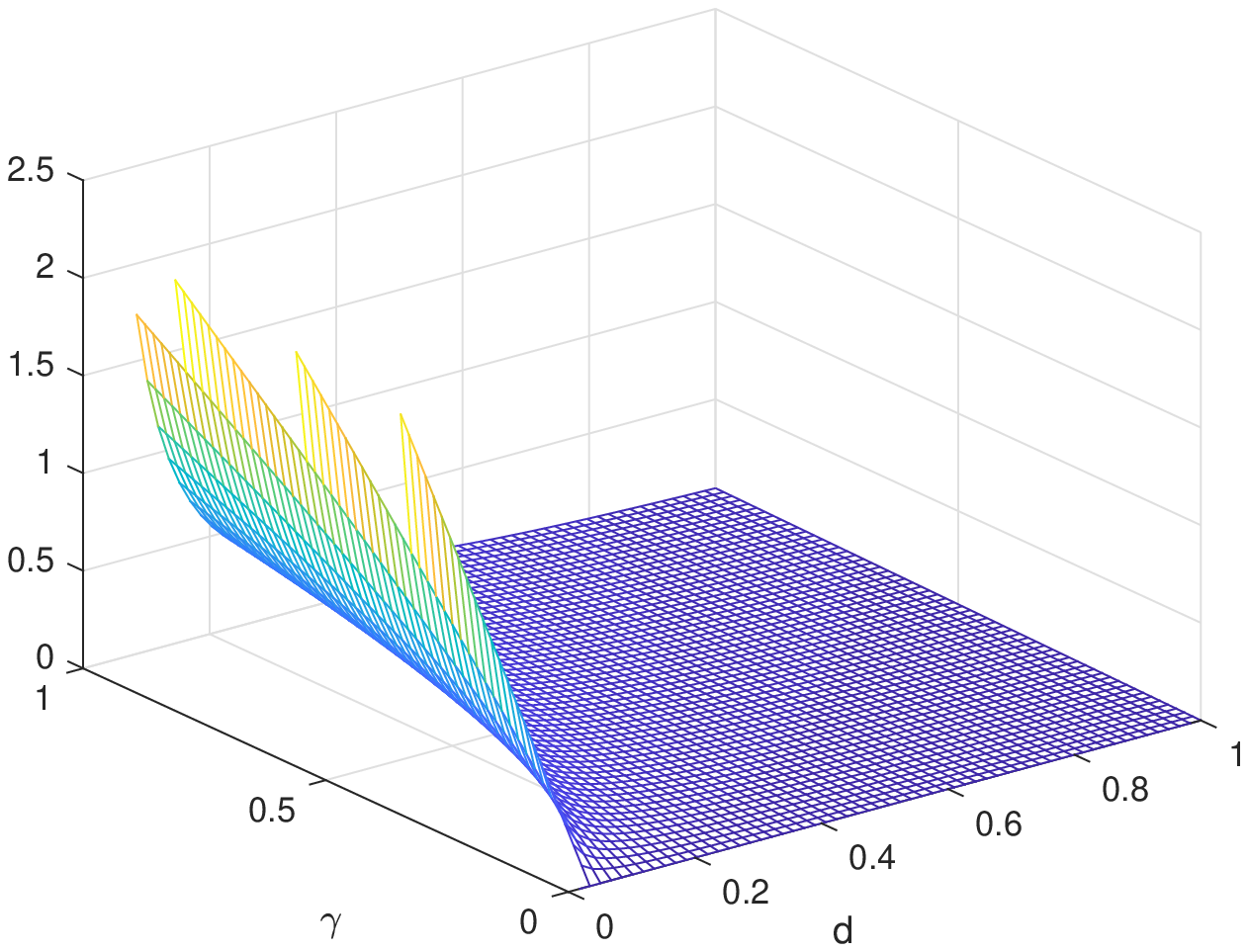}}
\subfloat[$R_{0}$ versus $\gamma$ and $\mu$]{
\label{fig9}\includegraphics[scale=0.45]{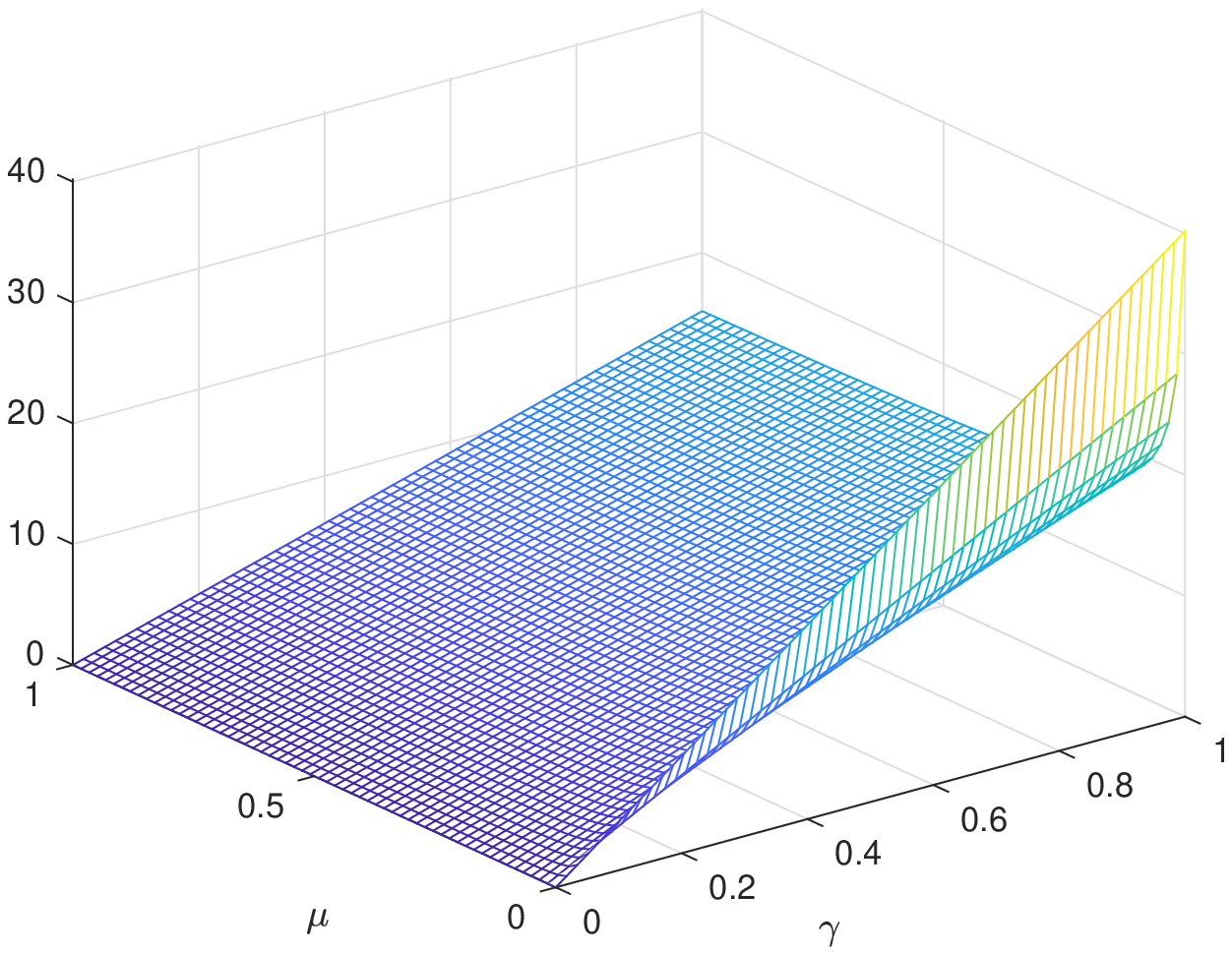}}\\
\subfloat[$R_{0}$ versus $\gamma$ and $\eta$]{
\label{fig10}\includegraphics[scale=0.45]{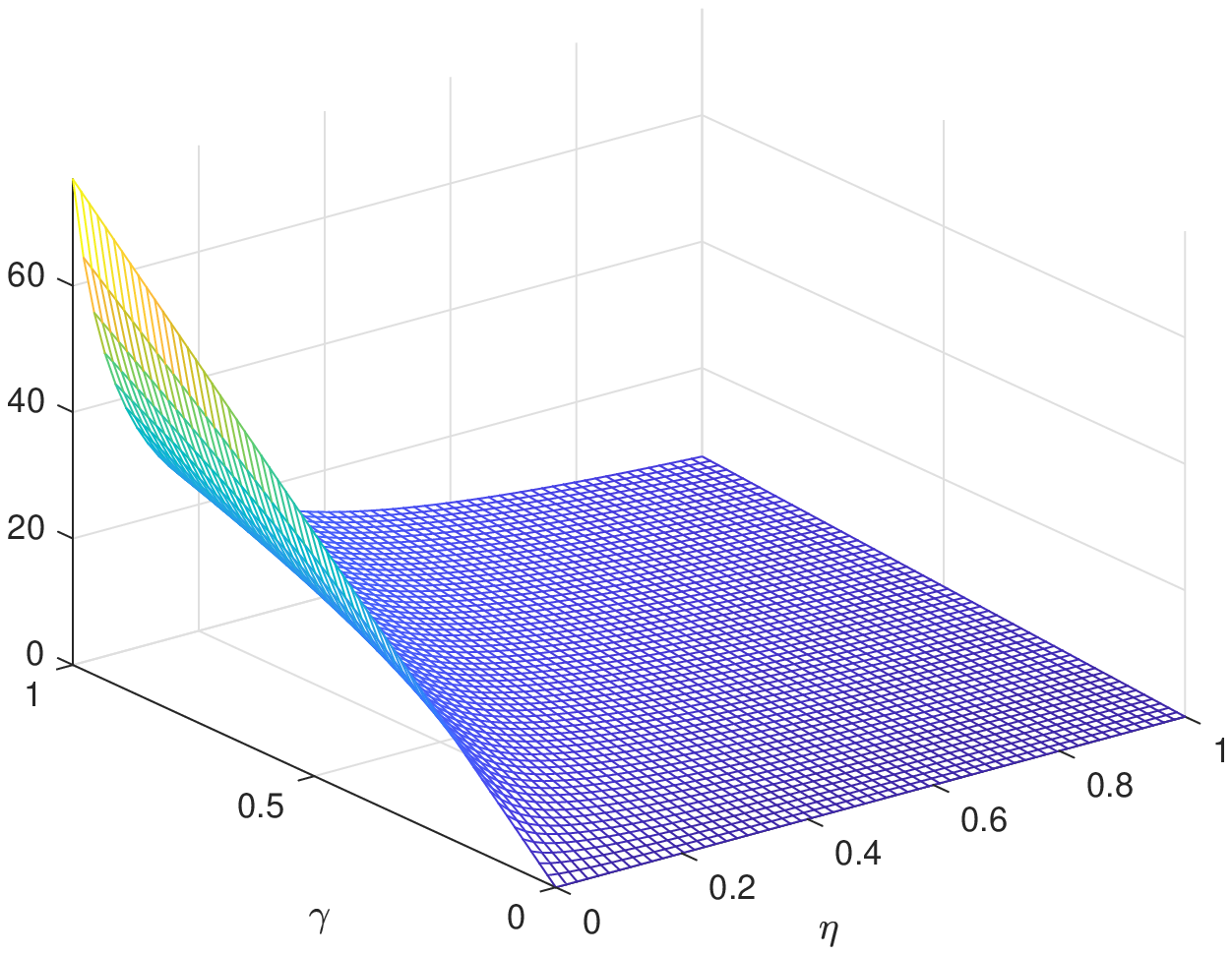}}
\subfloat[$R_{0}$ versus $h$ and $d$]{
\label{fig11}\includegraphics[scale=0.45]{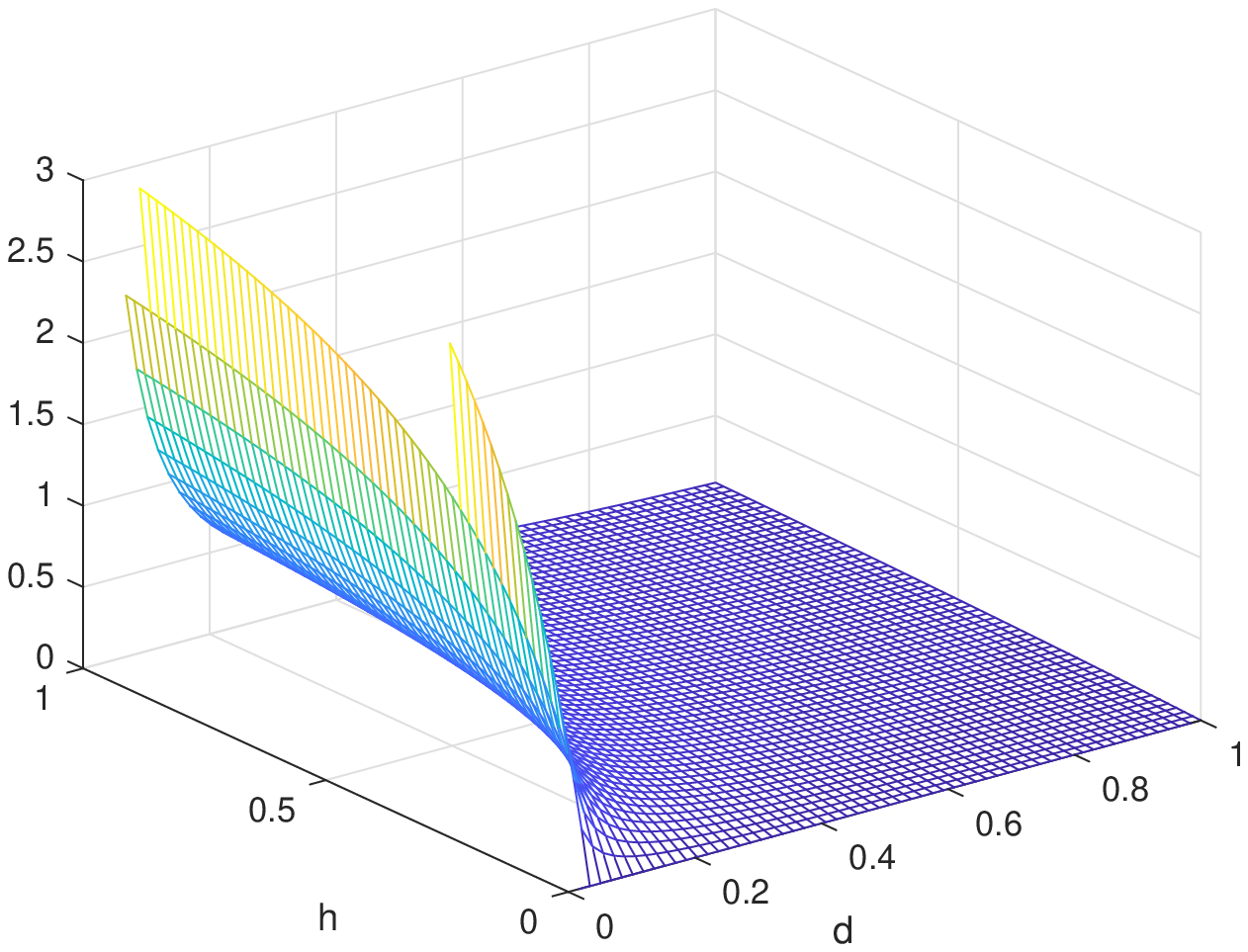}}\\
\subfloat[$R_{0}$ versus $h$ and $\mu$]{
\label{fig12}\includegraphics[scale=0.45]{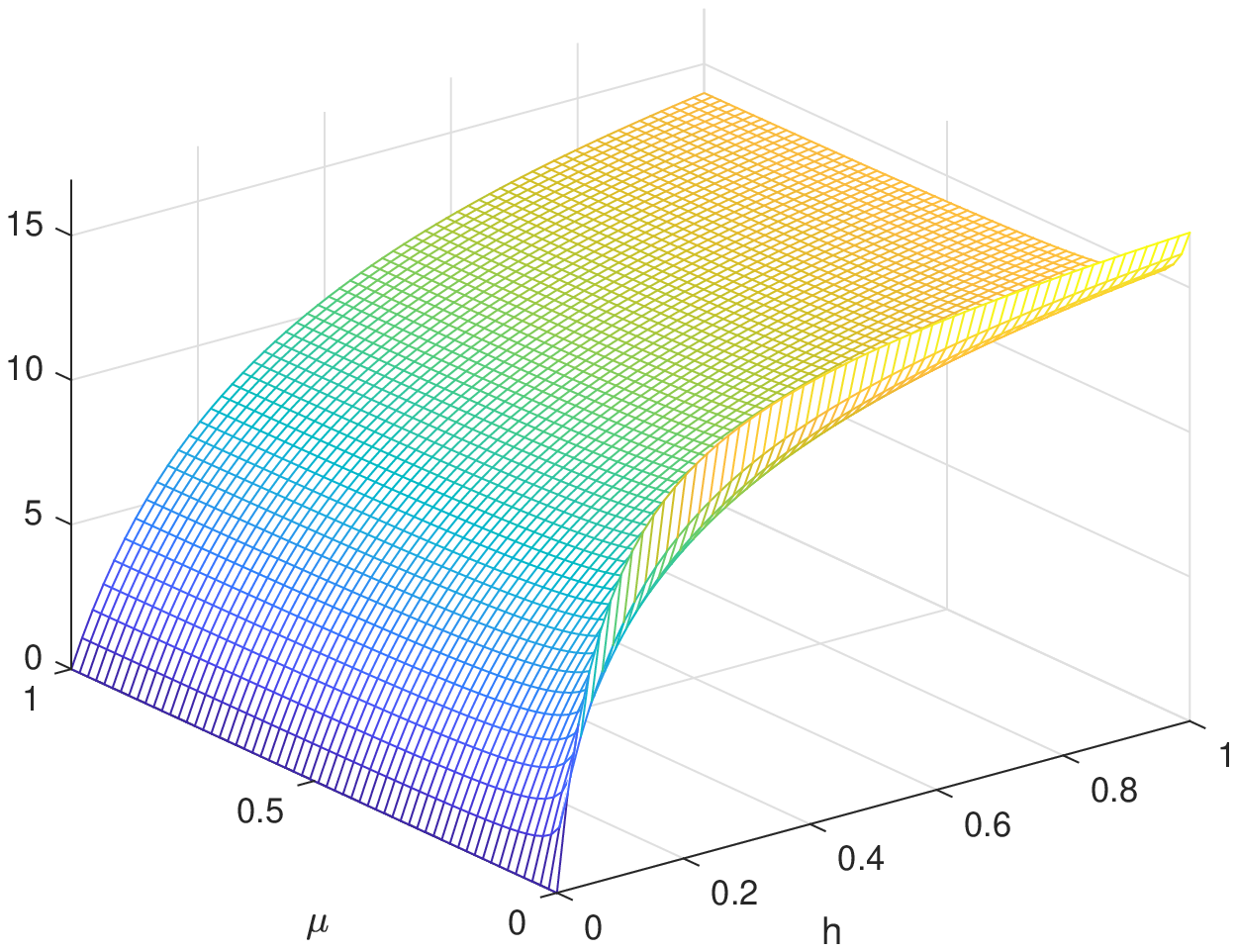}}
\subfloat[$R_{0}$ versus $h$ and $\eta$]{
\label{fig13}\includegraphics[scale=0.45]{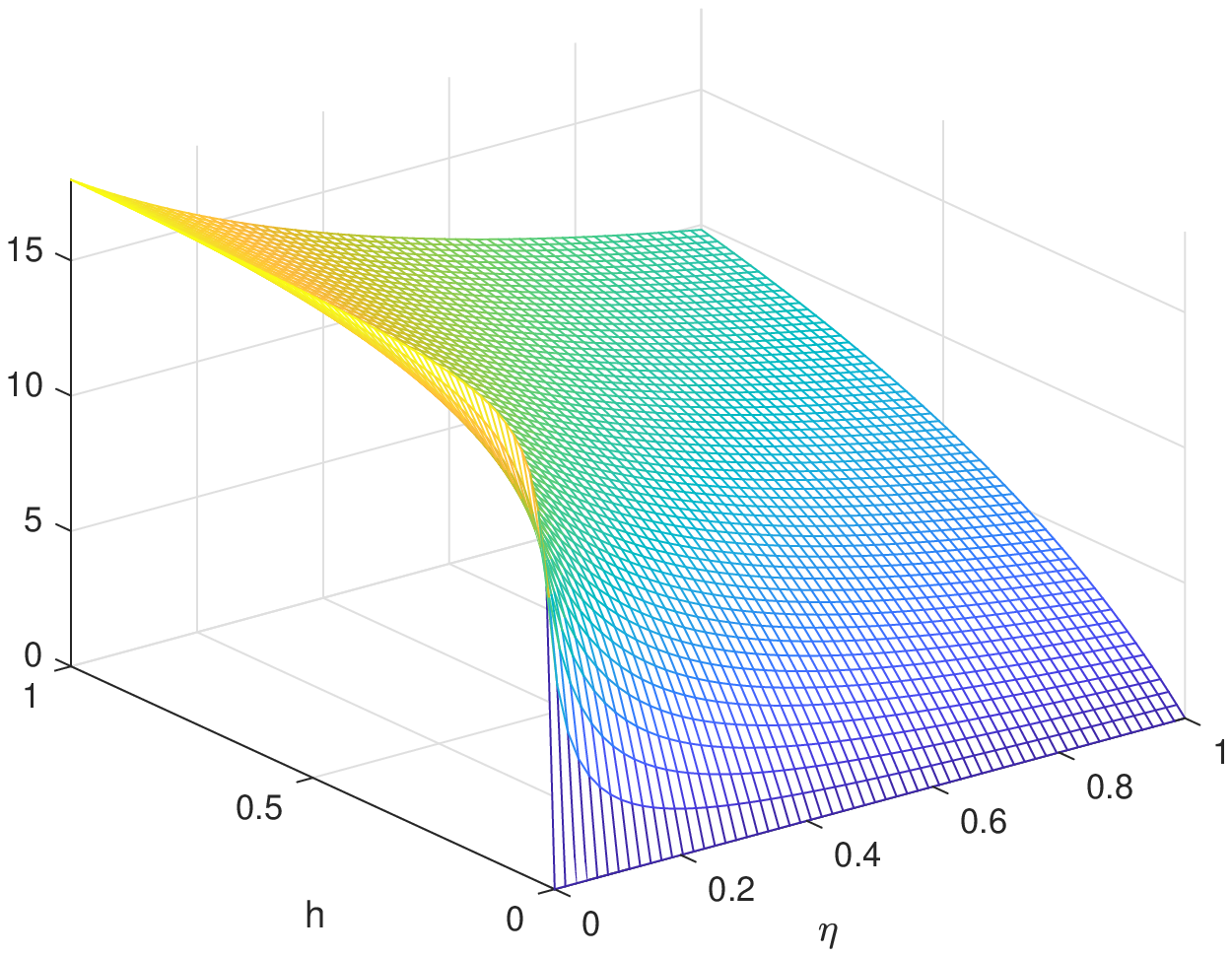}}\\
\subfloat[$R_{0}$ versus $d$ and $\sigma$]{
\label{fig15}\includegraphics[scale=0.45]{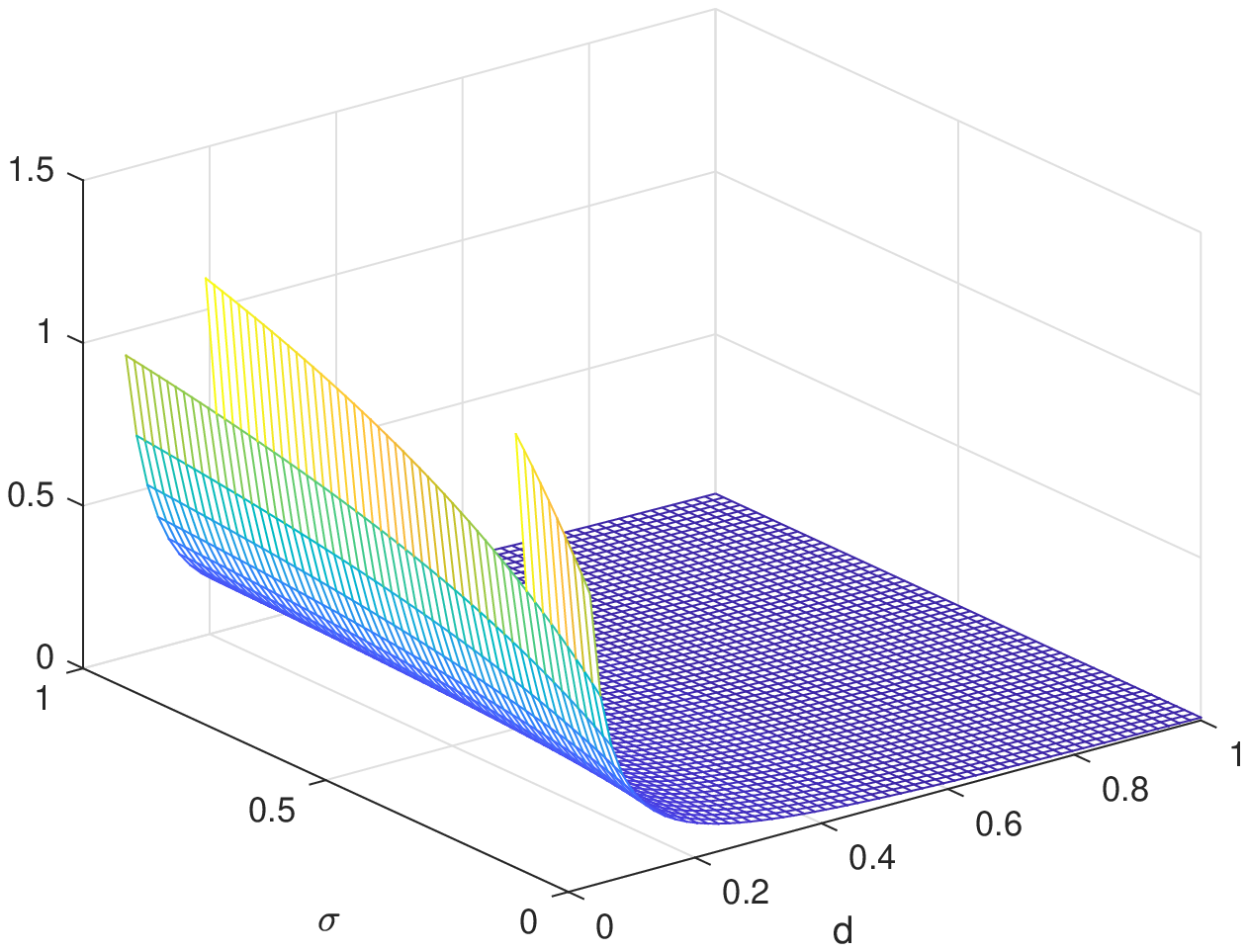}}
\subfloat[$R_{0}$ versus $d$ and $\eta$]{
\label{fig16}\includegraphics[scale=0.45]{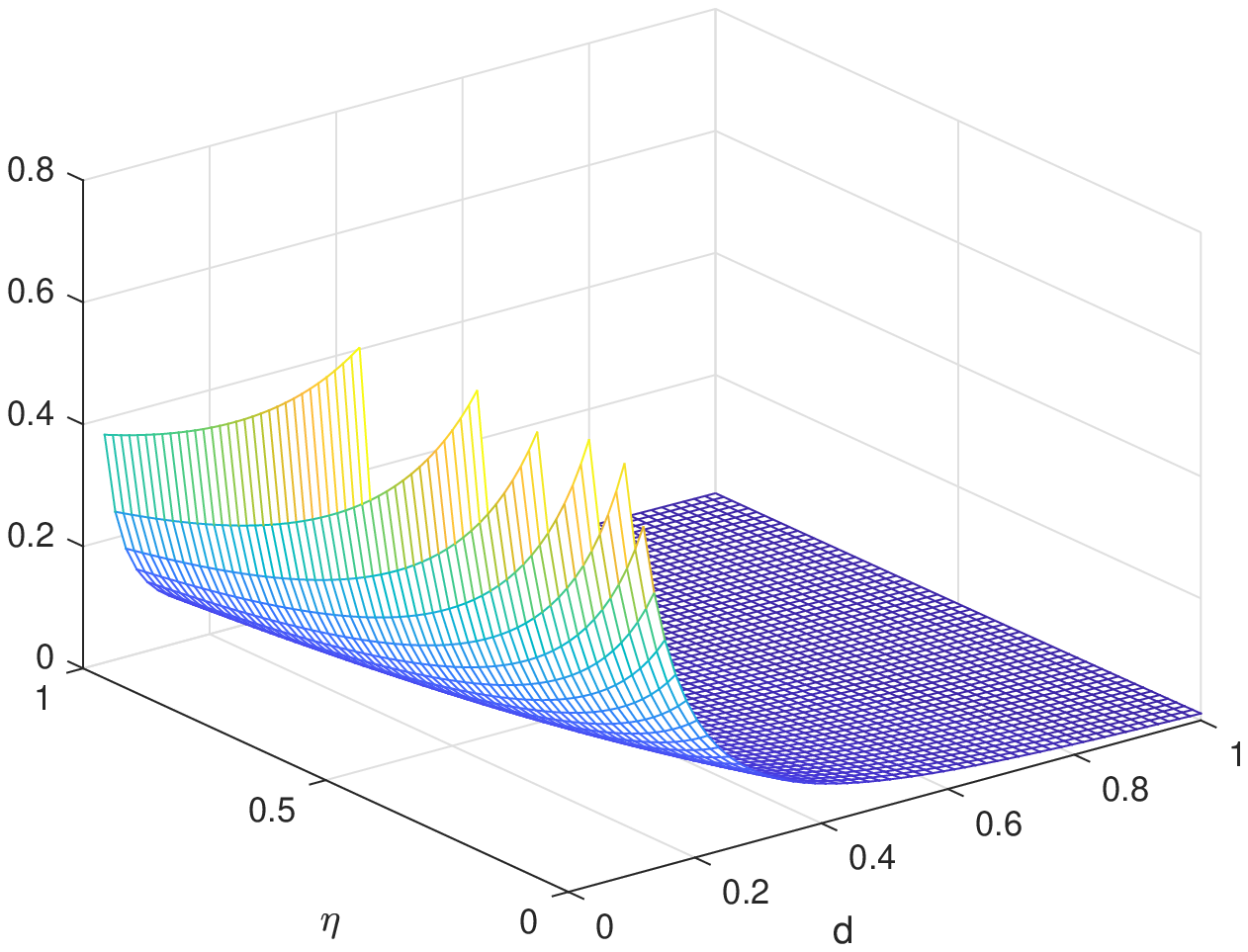}}}
\caption{Sensitivity of the basic reproduction 
number $R_0$ \eqref{eq2} for relevant parameters 
of model \eqref{m1}.}
\label{fig:sensitiv}
\end{figure}
\unskip


\section{Conclusions and Future~Work}
\label{sec:8}

In this manuscript, we studied a COVID-19 disease model
providing a detailed qualitative analysis and showed
its usefulness with a case study of Khyber Pakhtunkhawa, Pakistan.
Our sensitivity analysis shows that the
transmission rate $\gamma$ has a huge effect on the model as
compared to other parameters: the basic reproduction
number varies directly with the transmission rate $\gamma$. 
The sensitivity analysis also showed that 
the death rate parameter $\mu$ has no effect on
spreading the infection, which seems biologically correct. 
The~transmission rate will be small by keeping a social distancing 
and self-quarantine situation that causes a decrease in the infection.
In this way, one can control COVID-19 infection
from rapid spreading in the community. In~the future, we plan to
analyze optimal control techniques to reduce the population of
infected individuals by adopting a number of control measures.
A modification of the given model is also possible by introducing 
more parameters for analyzing the early outbreaks of COVID-19 and then
transmission and treatment aspects can be recalled. The~given system 
can be also simulated by adding exposed and hospitalized classes and
taking a stochastic fractional derivative.
Here, we have provided a case study 
with real data from Pakistan, but~other case studies can also be~done.

\vspace{6pt}
\authorcontributions{Conceptualization, M.R.S.A., 
A.K., A.Z., and D.F.M.T.; 
Formal analysis, M.R.S.A., A.D., A.K., and D.F.M.T.; 
Investigation, A.D. and A.Z.; Methodology, M.T. and A.Z.; 
Software, A.D. and A.K.; Supervision, D.F.M.T.; 
Validation, M.R.S.A. and D.F.M.T.; Writing---original draft, 
A.D., A.K., and A.Z.; Writing---review and editing, 
M.T., A.K., and D.F.M.T. 
All authors have read and agreed to the published version of the~manuscript.}

\funding{This research was partially funded by 
Funda\c{c}\~{a}o para a Ci\^{e}ncia e a Tecnologia (FCT) 
grant number UIDB/04106/2020 (CIDMA).}

\dataavailability{Not applicable.} 

\acknowledgments{The authors are grateful to reviewers for 
their comments, questions, and suggestions, which helped them
to improve the~manuscript.} 
	
\conflictsofinterest{The authors declare no conflicts of interest. 
The funders had no role in the design of the study; in the collection, analyses, 
or interpretation of data; in the writing of the manuscript; 
or in the decision to publish the~results.} 


\end{paracol}

\reftitle{References}


\end{document}